
\input amstex
\documentstyle{amsppt}
\magnification=\magstep1
\hcorrection{0in}
\vcorrection{0in}
\pagewidth{6.3truein}
\pageheight{9.3truein}
\NoRunningHeads

\def\Z{\Bbb Z}
\def\P{\Bbb P}
\def\O{\Cal O}
\def\Pic{\text{Pic }}
\topmatter
\title Higher dimensional examples of manifolds
whose adjoint bundles are not spanned
\endtitle
\author Takeshi Kawachi \endauthor
\affil Department of Mathematics, Tokyo Institute of Technology \endaffil
\address 2-12-1 Oh-okayama, Meguro-ku, Tokyo, JAPAN, 152 \endaddress
\email kawachi\@math.titech.ac.jp \endemail
\endtopmatter
\document
\head Introduction \endhead
Let $(X,L)$ be an $n$-dimensional polarized variety.
Fujita's conjecture says that if $L^n>1$ then the adjoint bundle
$K_X+nL$ is spanned and $K_X+(n+1)L$ is very ample.
There are some examples such that $K_X+nL$ is not spanned or
$K_X+(n+1)L$ is not very ample.
These are $(\P^n,\O(1))$,
hypersurface $M$ of degree $6$ in weighted projective space
$\P(3,2,1,1,\cdots ,1)$ with $\O_M(1)$ and numerically Godeaux surface etc.
Numerically Godeaux surface is the quotient space of a Fermat type hypersurface
of degree $5$ in $\P^3$ by an action of order $5$.
These examples are not so much.
We will construct new examples.\par
\proclaim{Theorem}
Let $n\geq 2$ be a natural number and $p\geq n+2$ be a prime number.
Let $X=(\xi_0^p+\cdots +\xi_{n+1}^p=0)$ be a Fermat type hypersurface
in degree $p$ in $\P^{n+1}$.
Let $\rho$ be a primitive $p$-th root of $1$.
Let $G$ be the cyclic group generated by $\rho$.
$G$ freely acts on $X$ by
$$ \rho\cdot (\xi_0:\cdots :\xi_{n+1})
  =(\rho^{k_0}\xi_0:\cdots :\rho^{k_{n+1}}\xi_{n+1}) $$
where $0\leq k_0<\cdots <k_{n+1}\leq p-1$.
Let $M$ be the quotient manifold of $X$ by $G$ and
$D$ be a divisor on $M$ whose degree is $1$.
If $\{k_0,\cdots,k_{n+1}\}$ satisfies
$k_i+k_j\equiv 2k_0+\sum_{t=0}^{n+1} k_t\mod{p}$
then $|K_M+nD|$ has a base point and $K_M+(n+1)D$ is not very ample.
\endproclaim\par
\head Introduction \endhead
Let $p\geq 5$ be a prime number and
$X$ be a Fermat type hypersurface of degree $p$ in $\P^{p-1}$.
If we take $(\xi_0:\cdots :\xi_{p-1})$ as coordinates of $\P^{p-1}$,
$X$ is defined by $\xi_0^p+\cdots +\xi_{p-1}^p=0$.
Let $\rho$ be a primitive $p$-th root of $1$
and $G$ be the cyclic group generated by $\rho$.
We define an action of $G$ on $\P^{p-1}$ by
$$\rho\cdot(\xi_0:\cdots:\xi_{p-1})=
  (\rho^0\xi_0:\cdots:\rho^{p-1}\xi_{p-1}).$$
Then this action is free, the quotient space $M=X/G$ is smooth
and the natural morphism $\pi\:X\to M$ is a covering of degree $p$.\par
\proclaim{Proposition}
The canonical sheaf of $M$ is trivial.
\endproclaim\par
\demo{Proof}
The canonical sheaf of $X$ is trivial and $\pi\: X\to M$ is unramified,
hence $\pi^*\omega_M\cong\omega_X=\O$.
If there exists $G$-invariant $p-2$ form on $X$,
that must be the pullback of a $p-2$ form on $M$.
Then we have a section of $\omega_M$ on $M$.
Hence we must show this.\par
Let $\{U_i\}$ be the affine open covering of $\P^{p-1}$ where
$U_i=\{\xi_i\ne 0\} \text{ for } i=0,\cdots,p-1$,
and $Z_i^j=\xi_j/\xi_i\text{ for } j\ne i$ be affine coordinates on $U_i$.
Since
$$X\cap U_i=\left\{1+\sum_{j\ne i} (Z_i^j)^p=0\right\}$$
on $U_i$, the case of $Z_i^0=\cdots =Z_i^{p-1}=0$ does not occur.
Hence there exists some $Z_i^j$ such that $Z_i^j\ne 0$.
Furthermore above equation shows that
$$\sum_{j\ne i} (Z_i^j)^{p-1} dZ_i^j=0.$$
Therefore $dZ_i^0,\cdots ,\check{dZ_i^j},\cdots ,dZ_i^{p-1}$ are
basis of $\Omega^1$.
Now we take $p-2$ form on $X\cap U_i\cap (Z_i^j\ne 0)$ as
$$\omega_{ij}=(-1)^{i+j+\delta(i<j)} \frac1{(Z_i^j)^{p-1}}
  dZ_i^0\wedge\cdots\wedge dZ_i^{p-1}
  \quad (\text{exclude } dZ_i^j) $$
where
$$\delta(\text{\it condition})=
  \cases 1, &  \text{if {\it condition} is true} \\
  0, & \text{if {\it condition} is false}.
  \endcases$$
Then this $\omega_{ij}$ is a $p-2$ form on
$X\cap U_i\cap (Z_i^j\ne 0)$ and does not equal to $0$ at anywhere.
On $X\cap U_i\cap (Z_i^j\ne 0)\cap (Z_i^k\ne 0)$,
we have
$$\align
  \omega_{ij} &= (-1)^{i+j+\delta(i<j)}\frac1{(Z_i^j)^{p-1}} dZ_i^0
    \wedge\cdots\wedge dZ_i^k\wedge\cdots\wedge dZ_i^{p-1} \\
  &= (-1)^{i+j+\delta(i<j)}\frac1{(Z_i^j)^{p-1}} dZ_i^0\wedge\cdots\wedge \\
    &\hskip3eM
    \left(\frac1{(Z_i^k)^{p-1}} (Z_i^0)^{p-1} dZ_i^0\wedge\cdots\wedge
    (Z_i^j)^{p-1} dZ_i^j\wedge\cdots\wedge
    (Z_i^{p-1})^{p-1} dZ_i^{p-1}\right) \\
    &\hskip3eM
    \wedge\cdots\wedge dZ_i^{p-1} \\
  &= (-1)^{i+j+\delta(i<j)}(-1)^{k-j+\delta(k<i<j\ or\ j<i<k)}
    \frac1{(Z_i^k)^{p-1}}dZ_i^0\wedge\cdots\wedge dZ_i^{p-1} \\
  &= (-1)^{i+k+\delta(i<j)+\delta(k<i<j\ or\ j<i<k)}
    \frac1{(Z_i^k)^{p-1}}dZ_i^0\wedge\cdots\wedge dZ_i^{p-1}. \\
\endalign$$
But $\delta$ satisfies
$$\delta(i<j)+\delta(k<i<j\ or\ j<i<k) =
  \cases 1, & \text{if $j<i<k$} \\
  1, & \text{if $i<k$ and (not $j<i$)} \\
  2, & \text{if $k<i<j$} \\
  0, & \text{if $k<i$ and (not $i<j$)}.
\endcases$$
Hence we have
$(-1)^{\delta(i<j)+\delta(k<i<j\ or\ j<i<k)}=(-1)^{\delta(i<k)}$.
Therefore we have
$$\omega_{ij}=(-1)^{i+k+\delta(i<k)}
  \frac1{(Z_i^k)^{p-1}}dZ_i^0\wedge\cdots\wedge dZ_i^{p-1}
  =\omega_{ik}.$$
Hence we can glue $\omega_{ij}$ in $U_i$,
we take a $p-2$ form $\omega_i$ on $U_i$
which does not equal to $0$ at anywhere on $U_i$.\par
On $X\cap U_i\cap U_j$, since
$$dZ_i^j=-(Z_j^i)^{-2} dZ_j^i,\quad
dZ_i^k=(Z_j^i)^{-1} dZ_j^k-Z_j^k (Z_j^i)^{-2} dZ_j^i$$
we have
$$\align
  \omega_{ij} &= (-1)^{i+j+\delta(i<j)}
    \frac1{(Z_i^j)^{p-1}} dZ_i^0
    \wedge\cdots\wedge\check{dZ_i^j}\wedge\cdots\wedge dZ_i^{p-1} \\
  &= (-1)^{i+j+\delta(i<j)} (Z_j^i)^{p-1} \bigwedge_{k\ne i,j}
    \left((Z_j^i)^{-1} dZ_j^k-Z_j^k (Z_j^i)^{-2} dZ_j^i \right) \\
  &= (-1)^{i+j+\delta(i<j)} (Z_j^i)^{p-1} \\
    &\hskip3eM \left(
    (Z_j^i)^{-(p-2)}\bigwedge_{k\ne i,j} dZ_j^k -
    (Z_j^i)^{-(p-1)}\sum_{\ell\ne j} Z_j^\ell
      \bigwedge_{k\ne\ell, j} dZ_j^k
    \right) \\
  &= (-1)^{i+j+\delta(i<j)} (Z_j^i)^{p-1} \\
    &\hskip3eM \left(
    (Z_j^i)^{-(p-2)}\bigwedge_{k\ne i,j} dZ_j^k +
    (Z_j^i)^{-(p-1)}\sum_{\ell\ne j} Z_j^\ell
      \frac{(Z_j^\ell)^{p-1}}{(Z_j^i)^{p-1}}\bigwedge_{k\ne i,j} dZ_j^k
    \right) \\
  &= (-1)^{i+j+\delta(i<j)} (Z_j^i)^{p-1} \\
    &\hskip3eM \left(
    (Z_j^i)^{-(p-2)}\bigwedge_{k\ne i,j} dZ_j^k +
    (Z_j^i)^{-2(p-1)}\sum_{\ell\ne j} (Z_j^\ell)^p\bigwedge_{k\ne i,j} dZ_j^k
    \right) \\
  &= (-1)^{i+j+\delta(i<j)} (Z_j^i)^{-(p-1)}\left(
    \sum_{\ell\ne j} (Z_j^\ell)^p \bigwedge_{k\ne i,j} dZ_j^k
    \right) \\
  &= -(-1)^{i+j+\delta(i<j)} (Z_j^i)^{-(p-1)} \bigwedge_{k\ne i,j} dZ_j^k \\
  &= (-1)^{i+j+\delta(j<i)} (Z_j^i)^{-(p-1)} \bigwedge_{k\ne i,j} dZ_j^k \\
  &= \omega_{ji}.
\endalign$$
Hence we can glue $\omega_{i}$, there exists a $p-2$ form $\omega$ on $X$
which does not equal to $0$ at anywhere on $X$.\par
Now we show this $\omega$ is $G$-invariant.
Since $X\cap U_i\cap (Z_i^j\ne 0)$ is $G$-invariant subspace of $X$
and $dZ_i^j$ is translated to $\rho^{j-i} dZ_i^j$ by $\rho\in G$
on this subspace.
Hence the action of $\rho$ translates $\omega_{ij}$ to
$$\align
  \rho\cdot\omega_{ij} &=
    (-1)^{i+j+\delta(i<j)}
    \rho^{-(p-1)(j-i)}(Z_i^j)^{-(p-1)}\bigwedge_{k\ne i,j}\rho^{k-i}dZ_i^k  \\
  &= (-1)^{i+j+\delta(i<j)}
     \rho^{\sum_{k\ne i}(k-i)}(Z_i^j)^{-(p-1)}\bigwedge_{k\ne i,j}dZ_i^k \\
  &= \omega_{ij}.
\endalign$$
Therefore $\omega$ is a $G$-invariant $p-2$ form.
Hence the canonical sheaf of $M$ is trivial.\qed
\enddemo\par
\head Fundamental case.\endhead
Let $H_j$ be a divisor on $X$ defined by $\xi_j=0$.
Then $G$ acts on $H_j$ freely.
Hence $D_j=H_j/G$ is a divisor on $M$.
Since $H_j^{p-2}=p$ and $\pi^*D_j=H_j$, we have $D_j^{p-1}=1$.\par
Let $N=D_1-D_0$ be a divisor on $M$.
Let $N'=D_j-D_{j-1}$.
Since $\pi^*(N'-N)=H_j-H_{j-1}-H_1+H_0$,
$\frac{\xi_{j-1}\xi_1}{\xi_j\xi_0}\in H^0(X,H_j-H_{j-1}-H_1+H_0)$
is a $G$ invariant rational function on $X$.
Hence $\frac{\xi_{j-1}\xi_1}{\xi_j\xi_0}$ is the pull back of a
rational function of $M$.
Since $\pi^*\O_M(N'-N)\cong\O_X$, we have $\O_M(N'-N)=\O_M$,
$N'$ and $N$ are linearly equivalent.
Hence we have $D_i-D_j\sim (i-j)N \mod{p}$.
In particular $N$ satisfies $\pi^*\O(N)\cong\O$ and $p\O(N)\cong\O$,
therefore $\O(N)$ is a torsion element of $\Pic(M)$.\par
Let $E=\sum_{i=0}^{p-1} a_i D_i+jN$ be a divisor on $M$
where each $a_i\geq 0$ and $j\geq 0$.
Since $D_i-D_0=iN$, we can write $E=tD_0+kN$
where $t=\sum_{i=0}^{p-1} a_i$ and $k\geq 0$.
Let $D$ be $D_0$, we consider divisors of form $tD+jN$.
\proclaim{Theorem 1}
Let $M$ be as above.
For $a\ne b$, Let $x_{a,b}$ be a point on $M$ defined by
$\xi_i=0$ for all $i\ne a,b$.
\roster
\item $|K_M+(p-2)D+jN|$ has finite base points.
  These base points are only $x_{a,b}$ where $a+b\equiv -j\mod{p}$.
\item $K_M+(p-1)D+jN$ is not very ample.
\endroster
\endproclaim\par
\demo{Proof}
Let $H$ be $H_0$.
(1) Since $\pi^*\: H^0((p-2)D+jN)\to H^0((p-2)H)$ is injective,
we identify $H^0((p-2)D+jN)$ to its image of $\pi^*$.
Then the basis of $H^0((p-2)D+jN)$ are forms $\xi_{i_1}\cdots\xi_{i_{p-2}}$
such that $i_1+\cdots+i_{p-2}\equiv j\mod{p}$.
Hence let $|\xi_{i_1}\cdots\xi_{i_{p-2}}|=i_1+\cdots+i_{p-2}$,
we define
$$S=\{f=\xi_{i_1}\cdots\xi_{i_{p-2}}\mid
  |f|\equiv j\mod{p}\}.$$
Let $a\ne b$.
We consider $x_{a,b}$ as one of $\pi^{-1}(x_{a,b})$,
let $x_{a,b}\in X$ be a point such that
$\xi_a=1,\xi_b=-1,\xi_i=0\text{ for }i\ne a,b$.
Let $a,b$ be satisfy $a+b\equiv -j$.
Suppose that $\pi(x_{a,b})$ is not a base point of $K_M+D+jN$.
Then there exists some $f\in S$ such that $f(x_{a,b})\ne 0$.
Since $x_{a,b}$ satisfies $\xi_i=0$ for $i\ne a,b$,
$f$ must be the form
$f_{a,b}^{(i)}=\xi_a^i\xi_b^{p-2-i}$ for $i=0,\cdots,p-2$.
Hence
$$|f_{a,b}^{(i)}|=ai+b(p-2-i)\equiv (a-b)i-2b \mod{p}.$$
Since $a\ne b$ and $0\leq a,b \leq p-1$, $a-b$ is prime to $p$.
Hence we have
$$\{(a-b)i-2b\mod{p}\mid i=0,\cdots ,p-1\} =\{0,\cdots ,p-1\}.$$
In particular, if $i=p-1$ then $(a-b)(p-1)-2b\equiv -a-b\equiv j\mod{p}$.
Therefore if $i=0,\cdots,p-2$ then
$|f_{a,b}^{(i)}|\equiv (a-b)i-2b\not\equiv j$,
hence $f_{a,b}^{(i)}\not\in S$.
This is contradiction to our assumption,
therefore points $x_{a,b}$ satisfied $a+b\equiv -j$ are
base points of $|K_M+(p-2)D+jN|$.\par
Let $q$ be $-j/2\mod{p}$ and
$x=(\xi_0,\cdots ,\xi_{p-1})$ is a base point of $|K_M+(p-2)D+jN|$.
Since the form $\xi_q^p-2$ is in $S$, $\xi_q$ must be $0$.
Now we assume $\xi_q=\cdots =\xi_{q+i}=0$.
Consider the form $f=\xi_{q+i+1}^a\xi_r^{p-2-a}$ where $r\ne q+i+1$.
Since $f$ in $S$, $a$ must be satisfy $(q+i+1-r)a=j+2r$.
If $r=-j-(q+i+1)$ then we have $a=p-1$, that is contradiction.
Hence if $\xi_{q+i+1}\ne 0$ then all other $\xi_r=0$ except $r=-j-(q+i+1)$.
If $\xi_{q+i+1}=0$ then we apply this argument as $i=i+1$.
This shows that base points of $|K_M+(p-2)D+jN|$ are
the points of form $\xi_{a,b}$ where $a+b\equiv -j\mod{p}$.
Hence the number of base points is $\frac{p-1}2$.\par
(2) Same as (1),
we identify the basis of $H^0((p-1)D+jN)$ are the form
$\xi_{i_1}\cdots\xi_{i_{p-1}}$ such that $i_1+\cdots+i_{p-1}\equiv j\mod{p}$.
Let $S$ be
$$S=\{f=\xi_{i_1}\cdots\xi_{i_{p-1}}\mid
  |f|\equiv j\mod{p}\}.$$
Forms in $S$ which is degree $1$ near $x_{a,b}$ must be the form
$\xi_c\xi_a^i\xi_b^{p-2-i}$ where $c\ne a,b,\ i=0,\cdots,p-2$.
Since
$$\align
  |\xi_c\xi_a^i\xi_b^{p-2-i}| &=c+ai+b(p-2-i) \\
  &\equiv c-2b+(a-b)i\equiv j\mod{p}. \\
\endalign$$
for any $i=0,\cdots,p-2$, we must be $c\equiv j+2b-(a-b)i \mod{p}$.
Since $a-b$ is prime to $p$, we have
$$\{j+2b-(a-b)i\mid i=0,\cdots,p-1\}=\{0,\cdots,p-1\}$$
and if $c\equiv j+2b-(a-b)(p-1)\equiv j+a+b$, then
$\xi_0\xi_a^i\xi_b^{p-2-i}\not\in S$ for any $i=0,\cdots,p-2$.
It said that
tangent $\frac{\partial}{\partial\xi_c}$ is not separated by
$|(p-1)D+jN|$ at $x_{a,b}$.
Therefore $K_M+(p-1)D+jN$ is not vary ample.\qed
\enddemo\par
We have following corollary if we take $j=0$. \par
\proclaim{Corollary}
$|K_M+(p-2)D|$ has base points and
$K_M+(p-1)D$ is not very ample.
\endproclaim\par
\head General case.\endhead
If we restrict above examples to a divisor $D_i$,
we can get examples of dimension $p-3$.
This process will give us examples of any dimension.\par
Let $n\geq 2$ be a natural number and let $p>n+2$ be a prime number.
Let $X'$ be a Fermat type hypersurface of degree $p$ in $\P^{n+1}$
defined by $\xi_0^p+\cdots+\xi_{n+1}^p=0$.
Let $G=<\rho>\cong \Z/p\Z$ be as above.
Let $0\leq k_0<\cdots <k_{n+1}\leq p-1$ be integers.
We define an action of $G$ on $\P^{n+1}$ by
$$\rho\cdot(\xi_0:\cdots:\xi_{n+1})=
  (\rho^{k_0}\xi_0:\cdots:\rho^{k_{n+1}}\xi_{n+1}).$$
Then this action is free, the quotient space $M'=X'/G$ is smooth
and the natural morphism $\pi\:X'\to M'$ is a covering of degree $p$.\par
Let $X,M$ be as in fundamental case and
let $S=\{0,\cdots ,p-1\}\setminus\{k_0,\cdots ,k_n\}$ and $s=\# S$.
We can identify $M'$ as hyperplane section of $M$ on $\cap_{i\in S} D_i$.
Let $D'=D_{k_0}|_{M'}$ and $N'=N|_{M'}$ be a divisor on $M'$.
Then we have
$$\align
  K_{M'}+tD'+jN'
  &= \left(K_M+(s+t)D_{k_0}+\left(\sum_{i\in S}i-k_0\right)N\right)\Big|_{M'}
\\
  &= \left(K_M+(s+t)D+\left(\sum_{i\in S} i+tk_0\right)N\right)\Big|_{M'}
\endalign$$
by canonical bundle formula.
Hence $|K_{M'}+tD'+jN'|$ has a base point if and only if
$|(K_M+(s+t)D+j'N)|_{M'}|$ has a base point on $\cap_{i\in S} D_i$
where $j'=\sum_{i\in S} i+tk_0$.
Similarly, $K_{M'}+tD'+jN'$ is not very ample if
$(K_M+(s+t)D+j'N)|_{M'}$ can not separate a tangent
$\frac{\partial}{\partial\xi_c}$ where $c\in S$ at $x$ on $\cap_{i\in S} D_i$.
\par
We can get following theorem from above arguments.\par
\proclaim{Theorem 2}
Let $n\geq 2$ and let $M'$ be as above.
For $a\ne b$, Let $x_{a,b}$ be a point on $M'$ defined by
$\xi_i=0$ for all $i\ne a,b$.
Let $k=\sum_{i\in S} i$.
Assume that $k_i+k_j\equiv 2k_0-k\mod{p}$ for some $i,j$.
\roster
\item $|K_{M'}+nD'|$ has a base point.
\item $K_{M'}+(n+1)D'$ is not very ample.
\endroster
\endproclaim\par
\demo{Proof}
If $p=n+2$ then this theorem is the same as fundamental case.
Let $p>n+2$ and $j'=\sum_{i\in S} (i-k_0)$.\par
(1) By fundamental case,
we have $|K_M+(p-2)D+(k-2k_0)N|$ has a base point
$x_{a,b}$ where $a+b\equiv 2k_0-k$.
Since we assume that $k_i+k_j\equiv 2k_0+k\mod{p}$,
we can take such $a,b$ as $a=k_i,b=k_j$.
Hence $x_{a,b}$ is on $\cap_{i\in S} D_i$,
this point is a base point of $|K_{M'}+nD'|$.\par
(2) By fundamental case,
we have $K_M+(p-1)D+(k-k_0)N$ does not separate the tangent
$\frac{\partial}{\partial\xi_c}$ at $x_{a,b}$ where $c\equiv a+b+(k-k_0)$.
Hence if we take $x_{a,b}$ as $a=k_i,b=k_j$ then
we have $a+b\equiv 2k_0-k$ and the tangent
$\frac{\partial}{\partial\xi_{k_0}}$ does not separate at $x_{a,b}$.
Therefore $K_{M'}+(n+1)D'$ is not very ample.\par
\enddemo\par
If $n=p-2$ then the canonical sheaf of $M$ is trivial.
Hence this case gives examples of Kodaira dimension $0$.
If $n<p-2$ then the canonical sheaf of $M$ is very ample.
Hence this case gives examples of general type.\par
\enddocument